\documentclass[12pt,reqno]{amsart}
\usepackage{xurl}
\usepackage{amsfonts,amsmath,amsthm}
\usepackage{amssymb,epsfig}
\usepackage{enumerate} 
\usepackage{appendix}
\usepackage{comment}

\usepackage{color} 
\definecolor{vert}{rgb}{0,0.6,0}

\usepackage[text={425pt,650pt},centering]{geometry}
\usepackage{graphicx}
\usepackage{epsfig}
\usepackage{tikz}
\usetikzlibrary{3d,calc,shapes}

\usepackage{caption}
\usepackage{color} 
\definecolor{vert}{rgb}{0,0.6,0}

\numberwithin{figure}{section}

\theoremstyle{plain}
\newtheorem{thm}{Theorem}[section]
\newtheorem{ass}{Assumption}
\newtheorem{defn}{Definition}

\newtheorem{lem}[thm]{Lemma}

\theoremstyle{remark}
\newtheorem{rem}{\bf{Remark}}
\numberwithin{equation}{section}

\newcommand{\N}{\mathbb{N}}

\newcommand{\R}{\mathbb{R}}

\newcommand{\Z}{\mathbb{Z}}

\newcommand{\cA}{\mathcal{A}}

\newcommand{\ol}{\overline}

\usepackage{booktabs,array,url}

\newenvironment{correction}{\begingroup\color{black}}{\endgroup}

\begin{document}
	
			\title[Controlling discrete semilinear wave equations toward flocks]
	{Controlling discrete semilinear wave equations toward flocks}

	\author[S. STRIKWERDA, H. V. TRAN, M.-B. TRAN]
	{Sarah Strikwerda, Hung Vinh Tran, Minh-Binh Tran}

	\thanks{
		The work of SS is partially supported by NSF grant DMS-2037851.
		The work of HVT is partially supported by NSF grant DMS-2348305.
		The work of MBT is partially supported by NSF Grants  CAREER DMS-2303146, DMS-2204795, DMS-2306379.
	}

	\address[S. Strikwerda]
	{
		Department of Mathematics, 
		Colorado State University, Weber Building, 841 Oval Drive, Fort Collins, Colorado 80523, USA}
	\email{s.l.strikwerda@colostate.edu}
	
	\address[H. V. Tran]
	{
		Department of Mathematics, 
		University of Wisconsin Madison, Van Vleck Hall, 480 Lincoln Drive, Madison, Wisconsin 53706, USA}
	\email{hung@math.wisc.edu}

	\address[M.-B. Tran]
	{
		Department of Mathematics, 
		Texas A\&M University, College Station, TX 77843, USA}
	\email{minhbinh@tamu.edu}


	\begin{abstract}
We construct pointwise bounded controls that steer a finite, periodic, fully actuated active harmonic lattice with Rayleigh--Helmholtz velocity forcing to flock states in finite time. In our setting, an initial Lyapunov-feedback phase brings the velocities and elastic forces into a prescribed small neighborhood of zero. Then, subsequent trajectory planning and cancellation of the known drift achieve exact arrival at a moving flock. For prescribed stationary or moving profiles satisfying a strict control-admissibility condition, a two-stage interpolation reaches the \begin{correction} target.\end{correction} We provide explicit sufficient time intervals for the stages following the initial feedback phase.    Finally, we formulate the associated minimum-time problem and its formal Hamilton--Jacobi equation, leaving the continuity of the value function open.
\end{abstract}

\maketitle\color{black}

	\section{Introduction}
	
	\subsection{Brief overview}

Active lattices combine elastic interactions with a local mechanism that supplies energy to the motion. In this work, we consider harmonic nearest-neighbor lattices coupling together with the Rayleigh--Helmholtz velocity force
\[
f(v)=(\alpha-\beta v^2)v,\qquad \alpha,\beta>0,
\]
which injects energy at small nonzero speeds and dissipates energy at large speeds. Rayleigh-type mechanisms arise in sustained vibrations \cite{Rayleigh1883,Rayleigh1945}, active Brownian particles \cite{SchweitzerEbelingTilch1998,ErdmannEbelingSchimanskyGeierSchweitzer2000,RomanczukBarEbelingLindnerSchimanskyGeier2012,BuradaLindner2012}, and self-propelled particle systems \cite{DorsognaChuangBertozziChayes2006}. The Toda--Rayleigh lattices studied in \cite{MakarovDelRioEbelingVelarde2001,ChetverikovEbelingVelarde2006} provide particularly close modeling precedents. Our system replaces the Toda interaction with a harmonic elastic force. It has no displacement-dependent on-site potential and is invariant under a common translation of all positions. We therefore refer to it as an active harmonic lattice.

The relationship between oscillator chains and flocking has been studied from several perspectives. Asymptotic flocking with Rayleigh friction is analyzed in \cite{HaHaKim2010}, while a one-dimensional Cucker--Smale system is derived as a singular limit of a damped oscillator chain in \cite{HaSlemrod2010}. Bonding forces are introduced in \cite{ParkKimHa2010} to ensure separation and spatial cohesion. Finite-time flocking for nonlinear Cucker--Smale models with Rayleigh-type friction and discrete $p$-Laplacian interactions is studied in \cite{KimLeePark2023}. Harmonic attraction and coherent motion also appear in the noisy active-particle model considered in \cite{Erdmann2005}, while delayed position coupling is investigated both theoretically and experimentally in \cite{Szwaykowska2016}.

Discrete wave control has also been studied as a spatial approximation of continuum control problems \cite{DZ,Zuazua2} (see \cite{Zuazua} for a broader perspective). A separate line of research concerns energy transport in anharmonic chains, originating with \cite{FPTU} and developed further in \cite{Campbell2005FPU,Gallavotti2008FPU,Lepri2003Thermal,Dhar2008Heat,Peierls,Spohn2,Spohn2014}. Control problems for the associated limiting transport models have recently been studied in \cite{HPPT1,HPPT2}, building on earlier work on the stochastic KdV-type equation \cite{HannaniRosenzweigStaffilaniTran2022}.

Control of collective motion is likewise well established. Optimal consensus control for the Cucker--Smale model is formulated in \cite{Bailo2018}, while constrained controls that steer systems toward flocks and mills are constructed in \cite{CKRT}. These studies motivate the  question studied in this work: can a pointwise bounded control steer an active harmonic lattice from arbitrary initial data exactly to a prescribed flock profile in finite time? In our setting, a ``flock'' is an exact configuration with a common constant velocity and a fixed spatial profile. 

The lattice consists of a fixed finite number of sites with periodic boundary conditions. Since a separate scalar control acts at each site,  the system is fully actuated in the acceleration variables. The main analytical difficulty is to dominate the active drift while respecting the pointwise control bound, first through a Lyapunov feedback and then through sufficiently slow interpolation of the spatial profile. After the first phase, the controls cancel the known drift along a planned trajectory. Since those controls  do not incorporate tracking-error feedback, exact steering requires exact knowledge of the model parameters and exact initialization at each switching time.
One of our main analytical contributions is a phase-by-phase construction for finite-time steering under a pointwise control bound on a finite, periodic, fully actuated active harmonic lattice. In particular, a strict target margin   is used quantitatively to absorb the acceleration and velocity-force terms.

\subsection{Settings and main goals}
	Let $D,n\in \N$ be given natural numbers, and let $h=1/D$ be the fixed grid size.
	Denote by 
	\[
	\Lambda  = \left\{\frac{j}{D}\,:\, j\in \Z \right\}^n=\left\{jh\,:\, j\in \Z \right\}^n
	\]
	the discrete grid in $\R^n$. 

	\smallskip

We consider the controlled harmonic lattice with Rayleigh--Helmholtz active friction
\begin{equation}\label{eq:wave}
\begin{cases}
w_{tt}-\Delta_Dw=(\alpha-\beta|w_t|^2)w_t+u(l,t)&\text{in }\Lambda\times(0,\infty),\\
w(\cdot,0)=w_0,\quad w_t(\cdot,0)=w_1&\text{on }\Lambda.
\end{cases}
\end{equation}
The nearest-neighbor force is elastic, whereas the velocity-dependent active force is nonconservative.  
The equilibria of the velocity force are $0$ and $\pm\sqrt{\alpha/\beta}$.

	Here, $\alpha,\beta>0$ are given constants, and $w_0, w_1$ are the given initial data, which are $\Z^n$-periodic. Let $e_k$ for $k \in \{1, \ldots, n\}$ be the standard basis vectors.
	For $l\in \Lambda$, we write $l' \sim l$ if $l'-l = \pm h e_k$ for some $k\in \{1,\ldots, n\}$.
	The corresponding discrete Laplacian is
	\[
	\Delta_D w(l,t)=\Delta w(l,t) = \sum_{l' \sim l} \frac{w(l',t) - w(l,t)}{h^2}.
	\]
	When there is no confusion, we write $\Delta_D w(l,t)=\Delta w(l,t)$ for simplicity.
	The function $u:\Lambda \times [0,\infty) \to \R$, which is $\Z^n$-periodic in $l$, is our control, and we typically assume that
	\[
	u \in \cA_M =\{f \in L^\infty(\Lambda \times [0,\infty))\,:\, f \text{ is $\Z^n$-periodic in $l$, }\|f\|_{L^\infty} \leq M\}
	\]
	where $M>0$ is a given positive constant (to be specified later).
	
	As we are in the spatial $\Z^n$-periodic setting, we are concerned with the spatial $\Z^n$-periodic solution $w$ to \eqref{eq:wave}.
	In particular, our problem can be phrased on the discrete torus 
	\[
	\overline\Lambda  = \left\{0,\frac{1}{D},\ldots, \frac{D-1}{D} \right\}^n,
	\]
	where we identify $k+j/D\equiv j/D$ for $k,j\in \Z$.
	If a function $f:\Lambda \to \R$ is $\Z^n$-periodic, we can think of $f$ as a function from $\overline \Lambda$ to $\R$ as well, and vice versa. {Throughout, $h$ is fixed and there are $N=D^n$ independently controlled sites.} 
	
	In this paper, we are interested in using the control $u$ to drive the solution $w$ of \eqref{eq:wave} to a flock where the group velocity is constant.
	Looking at \eqref{eq:wave}, it is clear that the nonlinear term vanishes if 
	\[
	w_t=\bar v = \pm\sqrt{\frac{\alpha}{\beta}} \qquad \text{ or } \qquad w_t=\bar v =0.
	\]
	Further, by using the ansatz $w(l,t)=\bar w(l) + \bar v t$ in \eqref{eq:wave}, we see that 
	\begin{equation}\label{eq:bar-w}
		|-\Delta \bar w(l)| = |u(l,t)| \leq M \quad \text{ for } l\in \overline \Lambda,    
	\end{equation}
	which gives a compatibility condition.
	
For $\bar v\in\{0,\pm\sqrt{\alpha/\beta}\}$, the compatibility condition $\|\Delta\bar w\|_\infty\leq M$ is necessary and sufficient to \emph{maintain} the exact flock: necessity follows by substitution, and sufficiency follows by taking $u_l=-\Delta\bar w(l)$. The reachability theorems for a prescribed profile below require the strict inequality. Maintenance at a saturated target and reachability of that target are separate assertions.

	\begin{defn}\label{def:flock} We say that $\varphi:\overline{\Lambda} \times \R \to \R$ is a moving flock if
		\[
		\varphi(l,t) = \bar \varphi(l) + \bar v t \quad \text{ for } (l,t) \in \overline{\Lambda} \times \R,
		\]
		where $\bar v = \pm\sqrt{\frac{\alpha}{\beta}}$. 
		
		We say that $\varphi:\overline{\Lambda} \times \R \to \R$ is a stationary flock if
		\[
		\varphi(l,t) = \bar \varphi(l)  \quad \text{ for } (l,t) \in \overline{\Lambda} \times \R.
		\]
	\end{defn}
	By symmetry, we only need to consider the moving flocks with $\bar v =\sqrt{\alpha/\beta}$.
	The main questions of interest in our paper are as follows.
	\begin{itemize}
		\item[(1)] Let $\bar v =\sqrt{\alpha/\beta}$ and $w_0, w_1$ be $\Z^n$-periodic initial data.
		Can we find a control $u\in \cA_M$ and an associated $\varphi_1$ such that $w=\varphi_1$ for all $t\geq T$, where $T>0$ is sufficiently large and $\varphi_1$ is a moving flock?
		
		\item[(2)] Let $\bar v =\sqrt{\alpha/\beta}$ and $w_0, w_1$ be $\Z^n$-periodic initial data.
		Take an arbitrary moving flock $\varphi_2(l,t)=\bar \varphi_2(l)+\bar v t$.
		Can we find a control $u\in \cA_M$ such that $w-\varphi_2$ is constant for all $t\geq T$ where $T>0$ is sufficiently large?
		
		\item[(3)] Let $\bar v =0$ and $w_0, w_1$ be $\Z^n$-periodic initial data.
		Take an arbitrary stationary flock $\varphi_3(l,t)=\bar \varphi_3(l)$.
		Can we find a control $u\in \cA_M$ such that $w-\varphi_3$ is constant for all $t\geq T$ where $T>0$ is sufficiently large?
	\end{itemize}
	
	\subsection{Main results}
	Write $x_l(t) = w(l,t)$ and $v_l(t) = w_t(l,t)$ for $l\in \Lambda$ and $t\geq 0$.
	Then, \eqref{eq:wave} becomes a system of ODEs, for $(l,t)\in \Lambda \times (0,\infty)$,
	\begin{equation}\label{eq:sys}
		\begin{cases}
			\dot x_l(t) = v_l(t),\\
			\displaystyle \dot v_l(t) = \sum_{l' \sim l} \frac{x_{l'}(t) - x_l(t)}{h^2} + (\alpha - \beta |v_l(t)|^2) v_l(t) + u_l(t),
		\end{cases}
	\end{equation}
	with the given initial condition
	\[
	\begin{cases}
		x_l(0)=w_0(l),\\
		v_l(0)=w_1(l).
	\end{cases}
	\]
	
The state space is $\R^{2D^n}$. Section~\ref{sec:connect} formulates the corresponding minimum-time problem and explains the qualifications needed for its Hamilton--Jacobi interpretation.

	\begin{ass}
		{We assume that the bound $M$ defining the admissible set $\cA_M$ satisfies} 
		\begin{equation}\label{ass-1}
			M>\frac{2}{3\sqrt{3}}\sqrt{\frac{\alpha^3}{\beta}}.
		\end{equation}
	\end{ass}
	
We write $F_*:=\frac{2}{3\sqrt3}\sqrt{\alpha^3/\beta}$. This is the height of the positive self-acceleration hump, $F_* =\max_{r\geq0} f(r)$. The condition $M>F_*$ is sufficient for our construction and sharp for strict dissipation by this method. Remark~\ref{rem:optimal} proves a reachability obstruction for a single site or {spatially homogeneous initial data under a common control}, with a prescribed positive target velocity.

	\begin{thm}\label{thm:main1}
		Assume \eqref{ass-1}.
		Let $w_0, w_1$ be $\Z^n$-periodic initial data.
		Then, there exist $u \in \cA_M$ and a moving flock 
        \[
		\varphi_1(l,t)=\bar{x}_{1l}+(t-T)\bar{v}=\bar{x}_{1l}-\bar{v}T +\bar v t \quad \text{ for }(l,t) \in \overline{\Lambda} \times \R,
		\] where $\bar{x}_{1l}\in\mathbb{R}$ are constants such that 
		 $w=\varphi_1$ for all $t\geq T$ where $T$ is taken to be sufficiently large. 
	\end{thm}

	\begin{thm}\label{thm:main2}
		Assume \eqref{ass-1}.
		Let $w_0, w_1$ be $\Z^n$-periodic initial data.
		Let $\varphi_2(l,t)=\bar \varphi_2(l)+\bar v t$ be an arbitrary flock, where $\bar \varphi_2$ satisfies
		\begin{equation}\label{eq:bar-varphi-strict}
			|-\Delta \bar \varphi_2(l)| < M \quad \text{ for } l\in \overline \Lambda.
		\end{equation}
		Then, there exists $u \in \cA_M$ such that the solution to \eqref{eq:sys} satisfies 
		\[
		\begin{cases}
			x_l(t)=\begin{correction}\bar \varphi_2(l)+\bar{v}t\end{correction},\\
			v_l(t)=\bar{v}=\sqrt{\alpha/\beta}
		\end{cases}
		\]
        for $l\in \overline \Lambda$ and {$t\geq T$}.
		In particular, we have $\begin{correction} w=\varphi_2\end{correction}$ for all $t\geq T$ where $T$ is taken to be sufficiently large.
	\end{thm}
	
	{A separate interpolation with zero terminal velocity gives the following result, whose proof is written out in Section~\ref{sec:thm2-3}.}
	\begin{thm}\label{thm:main3}
		Assume \eqref{ass-1}.
		Let $w_0, w_1$ be $\Z^n$-periodic initial data.
		Let $\varphi_3(l,t)=\bar \varphi_3(l)$ be an arbitrary stationary flock, where $\bar \varphi_3$ satisfies
		\begin{equation}\label{eq:bar-varphi-strict-n}
			|-\Delta \bar \varphi_3(l)| < M \quad \text{ for } l\in \overline \Lambda.    
		\end{equation}
		Then, there exists $u \in \cA_M$  such that the solution to \eqref{eq:sys} satisfies
		\[
		\begin{cases}
			x_l(t)=\begin{correction}\bar \varphi_3(l)\end{correction},\\
			v_l(t)=\bar{v}=0,
		\end{cases}
		\]
        for $l\in \overline \Lambda$ and {$t\geq T$}.
		In particular, we have $w=\begin{correction}\varphi_3\end{correction}$ for all $t\geq T$ where $T$ is taken to be sufficiently large.
	\end{thm}

The results provide explicit controls for finite-time steering, without asserting time optimality. Our construction assumes a finite periodic lattice with fixed mesh size and an independent scalar control at each site. Our trajectory-planning stages require exact knowledge of the model parameters and exact initialization at each switching time.

	\subsection*{Organization of the paper}
	The paper is organized as follows.
	In Section \ref{sec:thm1}, we prove Theorem \ref{thm:main1}.
	The proofs of Theorems \ref{thm:main2}--\ref{thm:main3} are given in Section \ref{sec:thm2-3}.
	{Section~\ref{sec:connect} gives the minimum-time formulation; Section~\ref{sec:numerics} illustrates the construction numerically.}
	
	\subsection*{Acknowledgement} \textnormal{ We would like to sincerely thank the Editor-in-Chief and the Handling Editor for their time and effort in taking care of our paper. We would also like to thank the referees for the reports.}

	 	\section{Proof of Theorem \ref{thm:main1}}\label{sec:thm1}
	
		Our strategy to prove Theorem \ref{thm:main1} is as follows.
	We start by using a control $u$ that drives the system toward zero velocity. 
	To do so, we first design $u$ so that the particles' velocities become sufficiently small after some finite time.
	Once the system is close to zero velocity, we can find a control in $\cA_M$ that brings it to exactly zero velocity, a stationary flock, after another finite amount of time. 
	Finally, we apply a third type of control to steer the system from a stationary flock to a moving flock with the desired group velocity $\bar v =\sqrt{\alpha/\beta}$.
	It is important to note that all the controls designed are explicit in terms of $(x_l(t))_{l\in \overline\Lambda}$ and $(v_l(t))_{l\in \overline\Lambda}$.
	
	\smallskip
	
	In the process of driving the system to zero velocity, we will study the Lyapunov function
	\[
	V(t)=\frac{1}{2} \sum_{l \in \overline \Lambda}\left[ |v_l(t)|^2 + \frac{1}{2}\sum_{l'\sim l} \frac{(x_{l'}(t)-x_l(t))^2}{h^2}\right].
	\]
	It is immediate that $V\geq 0$.
	We will now prove
	\[
	\dot{V}(t)=\sum_{l \in \overline \Lambda}\left[ (\alpha - \beta |v_l(t)|^2)|v_l(t)|^2 + (u_l(t),v_l(t)) \right]
	\]
	which shows that by choosing $u_l(t)$, we can control $\dot{V}$. 
	For clarity in our writing, we use the usual way to write the dot product $(a,b)=ab$ for $a,b\in \R$.
	
	We have
	\begin{align*}
		\dot{V}(t)&=\sum_{l \in \overline \Lambda} \left[(\dot{v}_l(t),v_l(t))+\frac{1}{2} \sum_{l' \sim l}\left(\frac{x_{l'}(t)-x_l(t)}{h^2}, v_{l'}(t)-v_l(t)\right)\right]\\
		&=\sum_{l \in \overline \Lambda} \Bigg[\left(\sum_{l'\sim l} \frac{x_{l'}(t)-x_l(t)}{h^2} + (\alpha -\beta |v_l(t)|^2)v_l(t)+u_l(t), v_l(t)\right)  \\
		&\qquad\qquad\qquad+\frac{1}{2} \left(\sum_{l' \sim l} \frac{x_{l'}(t)-x_l(t)}{h^2}, -v_l(t)\right)
		+ \frac{1}{2}  \sum_{l' \sim l} \left(\frac{x_{l'}(t)-x_l(t)}{h^2}, v_{l'}(t)\right)\Bigg]\\
		&=\sum_{l \in \overline \Lambda} \left((\alpha - \beta |v_l(t)|^2)v_l(t) + u_l(t), v_l(t)\right)\\
		&=\sum_{l \in \overline \Lambda}\left[ (\alpha - \beta|v_l(t)|^2)|v_l(t)|^2 + (u_l(t), v_l(t))\right].
	\end{align*}
	Let $\gamma > \max\left\{1, \frac{1}{M}\sqrt{\frac{\alpha^3}{\beta}}\right\}$ be fixed and define 
	\[
	a_1=\frac{1}{\gamma}\sqrt{\frac{\alpha}{\beta}} \qquad \text{ and } \qquad a_2=\gamma\sqrt{\frac{\alpha}{\beta}}.
	\]
	When $t \in [0,T_0]$ for some $T_0$ to be defined later, we will use the control $u_l$ based on $v_l$ in the following way:
	\begin{equation}\label{ui1}
		u_l(v_l)=\begin{cases} 0 &\text{if } |v_l|\geq 2a_2,\\
			-M\frac{v_l}{|v_l|}\left(2-\frac{|v_l|}{a_2}\right) &\text{if } |v_l|\in (a_2, 2a_2),\\
			-M\frac{v_l}{|v_l|} &\text{if } |v_l| \in [a_1, a_2],\\
			-M\frac{v_l}{a_1} &\text{if } |v_l| <a_1.
		\end{cases}
	\end{equation}
	Of course, $u_l$ is Lipschitz in $v_l$.
	See Remark \ref{rem:simple-u-l} for a discussion on this choice of control. 
	For each time $t$, we can split the indices $l$ into 3 sets:
	\[
	\begin{cases}
		I_2(t)=\{l \in \overline \Lambda\,:\,|v_l(t)|>a_2\}, \\
		I_1(t)=\{l \in \overline \Lambda\,:\,|v_l(t)| \in [a_1,a_2]\},\\
		I_0(t)=\{l \in \overline \Lambda \,:\,|v_l(t)|<a_1\}.
	\end{cases}
	\]
	For $l \in I_2(t)$, $u_lv_l=0$ if $|v_l| \geq 2a_2$ or 
	\[
	u_lv_l=-M\frac{v_l}{|v_l|}\left(2-\frac{|v_l|}{a_2}\right)v_l=-M|v_l|\left(2-\frac{|v_l|}{a_2}\right)<0
	\]
	if $|v_l|\in (a_2, 2a_2)$. 
	As $\gamma >1$,
	\[
	(\alpha - \beta |v_l(t)|^2)|v_l(t)|^2<(\alpha - \beta a_2^2)|v_l(t)|^2{=}(\alpha - \gamma^2\alpha)|v_l(t)|^2<0.
	\]
	Hence,
	\[
	(\alpha - \beta |v_l(t)|^2)|v_l(t)|^2 + u_lv_l < (\alpha - \gamma ^2\alpha) a_2^2.
	\]
	We define 
	\[
	C_2:=(\gamma^2-1)\alpha a_2^2 >0 .
	\]
	
	For $l \in I_1(t)$, we see $u_l v_l = -M|v_l|<0$. 
	Also, 
	\begin{equation} \label{I1} 
		(\alpha - \beta|v_l|^2)|v_l|^2 + u_lv_l = |v_l|(\alpha |v_l| - \beta |v_l|^3 - M).
	\end{equation} 
	Consider the function $f(x)=\alpha x - \beta x^3$. 
	It is clear that
	\begin{equation}\label{fmax}
	\max_{x\geq 0}f(x)= f\left(\sqrt{\frac{\alpha}{3\beta}}\right)= \frac{2}{3\sqrt{3}}\sqrt{\frac{\alpha ^3}{\beta}}<M
	\end{equation}
	by our assumption \eqref{ass-1}.
	Applying this to \eqref{I1}, we see 
	\begin{equation} \label{I1bound}
	(\alpha - \beta|v_l|^2)|v_l|^2 + u_lv_l \le |v_l|\left(\frac{2}{3\sqrt{3}}\sqrt{\frac{\alpha^3}{\beta}}-M\right) \le a_1\left(\frac{2}{3\sqrt{3}}\sqrt{\frac{\alpha^3}{\beta}}-M\right).
	\end{equation}
	We define
	\[
	C_1:=a_1\left(M-\frac{2}{3\sqrt{3}}\sqrt{\frac{\alpha^3}{\beta}}\right) >0.
	\]
	
	Finally, for $l \in I_0(t)$, $u_l v_l=-M\frac{|v_l|^2}{a_1}$. 
	Since $\gamma > \frac{1}{M} \sqrt{\frac{\alpha^3}{\beta}}$, we have $a_1<\frac{M}{\alpha}$. 
	Thus,
	\[(\alpha - \beta|v_l|^2)|v_l|^2 + u_l v_l=|v_l|^2\left(\alpha - \beta |v_l|^2 -\frac{M}{a_1}\right)\leq -\beta|v_l|^4.\]
	
	Therefore, we have shown that when $u_l$ is as described in \eqref{ui1}, 
	\begin{equation}
		\dot{V}(t)\leq \sum_{l \in I_2(t)} -C_2 + \sum_{l \in I_1(t)} -C_1 + \sum_{l \in I_0(t)} - \beta|v_l(t)|^4.
	\end{equation}
	Thus, $V(t)$ is nonincreasing and, in particular, $V(t)$ is bounded from above by $V(0)$. 
	Since both terms in $V$ are nonnegative, this means they must both be bounded from above. 
	Additionally, for $l\in \overline \Lambda$,
	\[\sum_{l' \sim l} \frac{(x_{l'}(t)-x_l(t))^2}{h^2} \leq 4 V(0),
	\]
	which implies that
	\begin{equation}\label{eq:bound-x-l}
		\left|\sum_{l' \sim l} \frac{x_{l'}(t)-x_l(t)}{h^2} \right| \leq C.
	\end{equation}
	Therefore, we have
	\[
	|\dot{v}_l(t)| \leq C (1+ \alpha+|u_l(t)|)\leq C \quad \text{or}\quad |\dot{v}_l(t)|\leq C(1+\beta C^2 +|u_l(t)|)\leq C.
	\]
	Here, $C=C(M,V(0),h,\alpha,\beta)>0$.

\begin{lem}\label{4.1}
		The system \eqref{eq:sys} satisfies 
		\[
		\lim_{t \to \infty} v_l(t)=0 \quad \text{ for all $l\in \overline \Lambda$}
		\]
		when $u_l$ is given by \eqref{ui1}.
	\end{lem}
	\begin{proof}
		Assume this is not the case. 
		Then there exists a $c_0>0$, an index $l$, and a sequence of times $t_k \to \infty$ such that 
		\[
		|v_l(t_k)|>c_0.
		\]
		Since $|\dot{v}_l(t)| \leq C$, there is a $\tau>0$ such that 
		\[
		|v_l(t)|>\frac{c_0}{2} \quad \text{ for all $t \in (t_k-\tau, t_k+\tau)$, $k\in \N$.}
		\]
       Because $t_k\to\infty$, we may pass to a subsequence, still denoted by $t_k$, such that $t_1>\tau$ and $t_{k+1}-t_k>2\tau$ for every $k\geq1$. The intervals $(t_k-\tau,t_k+\tau)$ are therefore pairwise disjoint and contained in $(0,\infty)$. We then set
\[
c_*:=\min\{C_2,C_1,\beta c_0^4/16\}>0.
\]
        On each of these intervals, $\dot V\leq-c_*$, while $\dot V\leq0$ elsewhere. Therefore, for every $m\geq1$, we obtain
\[
0\leq V(t_m+\tau)
=V(0)+\int_0^{t_m+\tau}\dot V(s)\,ds
\leq V(0)-2m\tau c_*.
\]
As $m\to\infty$, the right-hand side tends to $-\infty$, contradicting the nonnegativity of $V$.
Thus, $\lim_{t\to\infty}v_l(t)=0$.
	\end{proof}
	
	\begin{rem}\label{rem:simple-u-l}
		In \eqref{ui1} above, we define $u_l=u_l(v_l)$ in such a way that we do not need to use the control $u_l$ when $|v_l|\geq 2a_2$. This choice is not necessary, but for practical purposes, having the control set to zero might be preferred, and allowing the system to be driven by the large velocities in the system might be preferable. Another way, which is simpler but requires a higher magnitude of control, is to define $u_l$ as
		\begin{equation*}
			u_l(v_l)=\begin{cases} 
				-M\frac{v_l}{|v_l|} &\text{if } |v_l|  \geq a_1,\\
				-M\frac{v_l}{a_1} &\text{if } |v_l| <a_1.
			\end{cases}
		\end{equation*}
		The proof that this control drives the system towards zero velocity would be identical to the proof above except for the definition of $I_1(t)$.
        
	\end{rem}

\begin{rem}\label{rem:optimal}
The condition $M>F_*$ ensures strict dissipation for the feedback above and also has a precise reachability interpretation in the scalar case. Let $r_*:=\sqrt{\alpha/(3\beta)}$, so that $f(-r_*)=-F_*$. We now consider the scalar equation
\[
\dot v=f(v)+u(t),\qquad |u(t)|\leq M,
\]
and assume that $M\leq F_*$. We then have 
\[
\dot v\le f(v)+M,
\]
which means $v$ is a subsolution.

Since $f(-r_*)+M\leq0$, the constant function $q(t)\equiv-r_*$ is a supersolution of $\dot q=f(q)+M$. 
Assume $v(0)\leq-r_*$.
We then deduce, by the comparison principle,
\[
 v(t)\leq-r_*
\qquad(t\geq0).
\]
Indeed, on any bounded time interval, we choose a Lipschitz constant $L_f$ for $f$ on an interval containing the range of $v$ and the point $-r_*$. Then, we can bound
\[
\frac{d}{dt}(v+r_*)_+\leq L_f(v+r_*)_+
\qquad\text{a.e.},
\]
Therefore, we find $(v(0)+r_*)_+=0$ and hence $(v(t)+r_*)_+=0$ throughout that interval.

Thus, initial velocities satisfying $v(0)\leq-r_*$ cannot reach the prescribed positive velocity $\sqrt{\alpha/\beta}$, including in the borderline case $M=F_*$. The same obstruction holds for spatially homogeneous initial positions and velocities when a common control is applied at every site. 
\end{rem}
	
	\begin{lem}
		The system \eqref{eq:sys} satisfies 
		\begin{equation}\label{eq:lim-x-l}
			\lim_{t \to \infty} \sum_{l' \sim l} \frac{x_{l'}(t)-x_l(t)}{h^2} =0 \quad \text{ for all $l \in \overline \Lambda$}
		\end{equation}
		when $u_l$ is given by \eqref{ui1}.
	\end{lem}
	
	\begin{proof}
		Since $|\dot x_l(t)| = |v_l(t)|\le \sqrt{2V(0)}$, the functions $x_l(t)$ are Lipschitz with constant $L:=\sqrt{2V(0)}$. 
		Then $x_{l'}(t)-x_l(t)$ are Lipschitz with Lipschitz constants less than or equal to $2L$. {If $V(0)=0$, all velocities and nearest-neighbor differences vanish, and the conclusion is straightforward. Hence we may assume $L>0$.} 
		Assume for the sake of contradiction that \eqref{eq:lim-x-l} does not hold.
		In light of \eqref{eq:bound-x-l}, there exists an $l$ and a sequence $\{t_k\}$ such that $t_k \to \infty$ and 
		\begin{equation} \label{limc}
		\lim_{ k \to \infty} \sum_{l' \sim l} \frac{x_{l'}(t_k)-x_l(t_k)}{h^2}=c \neq 0.
		\end{equation}
		There exists a constant $K>0$ such that when $k>K$
		\begin{align*}
			&\left|\sum_{l' \sim l} \frac{x_{l'}(t)-x_l(t)}{h^2}-c\right|\\
			=\, &\left|\left[\sum_{l' \sim l}\frac{(x_{l'}(t)-x_l(t))-(x_{l'}(t_k)-x_l(t_k))+(x_{l'}(t_k)-x_l(t_k))}{h^2}\right]-c\right|\\
			\leq\, & \frac{8n L}{h^2}|t_k-t|
            +\frac{|c|}{8}.\end{align*}
		Additionally, there exists a sequence $(\eta_k)_{k \geq K}$ such that when $t>t_k$, 
		\begin{equation}\label{2.9}
		\begin{cases}
			|v_l(t)|<\eta_k<a_1, \\
			|(\alpha-\beta|v_l(t)|^2)v_l(t)|<\eta_k<a_1
		\end{cases}
		\end{equation}
		and $\eta_k \to 0$ due to Lemma \ref{4.1}. 
		Thus,
		\begin{equation}
			\begin{split} \label{31}
				&|\dot v_l(t)-c|
				=\, \left|\sum_{l' \sim l} \frac{x_{l'}(t)-x_l(t)}{h^2}-c+(\alpha - \beta |v_l(t)|^2)v_l(t) + u_l(t)\right|\\
				\leq \, & \frac{8nL}{h^2}|t_k-t|
                +\frac{|c|}{8}+\eta_k +\frac{M}{a_1}\eta_k.
		\end{split} \end{equation}
		Let $\tau=\min\{1, \frac{|c|h^2}{16nL}\}$. 
		Using \eqref{31}, we obtain
		\begin{align*}
			&|v_l(t_k +\tau) - c\tau| - |v_l(t_k)| 
			\leq  |v_l(t_k + \tau) -v_l(t_k) - c\tau|\\
			=\, &\left|\int^{t_k+\tau}_{t_k} \left(\dot v_l(s) -c\right)\, ds\right|
			\leq\, \int^{t_k+\tau}_{t_k} |\dot v_l(s)-c|\,ds\\\leq &\tau\left(\frac{8nL\tau}{h^2}
            +\frac{|c|}{8}+\eta_k + \frac{M}{a_1}\eta_k\right)\leq\tau\left(\frac{|c|}{2}
            +\frac{|c|}{8}+\eta_k + \frac{M}{a_1}\eta_k\right) .
		\end{align*}
		Thus, 
		\begin{align*}
			&|v_l(t_k+\tau)-c\tau|
			\leq\, |v_l(t_k)|+\tau\left(\frac{|c|}{2} 
            +\frac{|c|}{8}+\eta_k +\frac{M}{a_1}\eta_k\right)\\
			\leq\,& \eta_k +\tau\left(\frac{|c|}{2} 
            +\frac{|c|}{8}+\eta_k +\frac{M}{a_1}\eta_k\right).
		\end{align*}
		Using the reverse triangle inequality again, we have
		\[|c \tau|-|v_l(t_k+\tau)| \leq \eta_k +\frac{|c|\tau}{2} 
        +\frac{|c|\tau}{8}+\left(1+\frac{M}{a_1}\right)\tau\eta_k,\]
        so
		\begin{equation} \label{167} 
			\frac{3|c| \tau}{8} -\tau \eta_k\left(1+\frac{M}{a_1}\right)-\eta_k\leq |v_l(t_k + \tau)| \quad \text{ for all } k >K.
		\end{equation}
		Since  $\eta_k \to 0$, we know there exists $K_2$ such that when $k>K_2$, $\eta_k(1+\frac{M}{a_1})\leq \frac{|c|}{16}$ and $\eta_k < \frac{|c|\tau}{16}$. Thus, when $k>\max\{K, K_2\}$, 
		\[
		|v_l(t_k + \tau)|\geq \frac{|c| \tau}{4}.
		\]
		However, if we take the limit as $k \to \infty$, we have $\frac{|c|\tau}{4}\leq 0$. 
		Therefore, $c=0$, which contradicts \eqref{limc}.
		The proof is complete.
	\end{proof}
	
	Therefore, there exists a $T_0$ such that if \eqref{ui1} is used for $0<t\leq T_0$, 
	\[
	|v_l(T_0)|<\epsilon \quad \text{ and } \quad \left|\sum_{l' \sim l}\frac{x_{l'}(T_0)-x_l(T_0)}{h^2}\right|<\epsilon
	\]
	for some $\epsilon$ satisfying 
	\begin{equation} \label{eps}
	0<\epsilon<\min\left\{1, \frac{M}{2n+1+\frac{8n}{h^2}+2\alpha+8\beta}, \frac{M}{2+\frac{16n}{h^2}}\right\}.
	\end{equation}
	
	It is important to note that $\epsilon$ can be chosen explicitly here with the above bound, which is useful in practice.
	We are now ready to finish the proof of Theorem \ref{thm:main1}.
	\begin{proof}[Proof of Theorem \ref{thm:main1}]
		Define $T_{1,l}=T_0 + \frac{1}{\epsilon}|v_l(T_0)|$ for each $l \in \overline \Lambda$ and $T_1=\max_{l\in \overline\Lambda}T_{1,l}$. 
		For $T_0<t\leq T_1$, we define 
		\begin{equation} \label{u2}u_l(t)=\dot{v}_l(t) - \sum_{l'\sim l} \frac{x_{l'}(t)-x_l(t)}{h^2} - (\alpha - \beta|v_l(t)|^2)v_l(t)\end{equation}
		where $v_l(t)$ satisfies the following:
		\begin{equation} 
			v_l(t)=\begin{cases}v_l(T_0)-\epsilon(t-T_0)\frac{v_l(T_0)}{|v_l(T_0)|} \quad &\text{for } t \in [T_0, T_{1,l}]{\text{  and  } v_l(T_0)\neq 0},\\
				0 &\text{for } t>T_{1,l}{\text{  or  } v_l(T_0)=0}.\end{cases}
		\end{equation}
		
	For \(v_l(T_0)\ne0\),
		\begin{align*}
			x_l(t)&=x_l(T_0)+\int_{T_0}^t \dot{x}_l(s) \,ds=x_l(T_0)+\int_{T_0}^t \left(v_l(T_0)-\epsilon(s-T_0)\frac{v_l(T_0)}{|v_l(T_0)|}\right)\, ds\\
			&=x_l(T_0)+v_l(T_0)(t-T_0)-\frac{\epsilon}{2} \frac{v_l(T_0)}{|v_l(T_0)|}(t-T_0)^2
		\end{align*}
		when $t \in [T_0, T_{1,l}]$ and $x_l(t)=x_l(T_{1,l})$ when $t>T_{1,l}$.
		Notice
		\begin{align*}
			|x_l(t)-x_l(T_0)|&\leq |v_l(T_0)|(T_{1,l}-T_0)+\frac{\epsilon}{2}(T_{1,l}-T_0)^2\\
			&\leq \frac{|v_l(T_0)|^2}{\epsilon} + \frac{\epsilon}{2}\frac{|v_l(T_0)|^2}{\epsilon^2}<2\epsilon.
		\end{align*}
		We must show that $|u_l|<M$.
		We first observe that when $t \in (T_{1,l},T_1]$,
		{\begin{align*}
			|u_l(t)|&=\left|\sum_{l' \sim l} \frac{x_{l'}(t)-x_l(t)}{h^2}\right|\\
			&=\left|\sum_{l' \sim l} \frac{x_{l'}(t)-x_{l'}(T_0)+x_{l'}(T_0)-x_l(T_0)+x_l(T_0)-x_l(t)}{h^2}\right|\\
			&\leq \sum_{l'\sim l}\left(\frac{2\epsilon}{h^2} +\frac{2\epsilon}{h^2}\right)+\left|\sum_{l'\sim l}\frac{x_{l'}(T_0)-x_l(T_0)}{h^2}\right|\leq \frac{8n\epsilon}{h^2} +\epsilon
		\end{align*}}
		when $t \in (T_{1,l}, T_1].$ Similarly, when $t\in (T_0, T_{1,l}]$,
		\begin{align*} 
			u_l(t)&=-\epsilon\frac{v_l(T_0)}{|v_l(T_0)|}-\sum_{l' \sim l} \frac{x_{l'}(t)-x_l(t)}{h^2} \\
			&\qquad-\left(\alpha - \beta\left|v_l(T_0)-\epsilon(t-T_0)\frac{v_l(T_0)}{|v_l(T_0)|}\right|^2\right)\left(v_l(T_0)-\epsilon(t-T_0)\frac{v_l(T_0)}{|v_l(T_0)|}\right).
		\end{align*}
		Hence, using the bound on $\epsilon$ given in \eqref{eps} we see
		\begin{align*}
			&|u_l(t)|\\
			\leq\, &\epsilon + \left|\sum_{l' \sim l} \frac{x_{l'}(t)-x_l(t)}{h^2} \right| +\alpha\left|v_l(T_0)\right|+\alpha \left|\epsilon(T_{1,l}-T_0)\right| + \beta \left(\left|v_l(T_0)\right| +\epsilon|T_{1,l}-T_0|\right)^3\\
			\leq \, &\epsilon + \frac{8n\epsilon}{h^2}+2n\epsilon +\alpha\epsilon+\alpha \epsilon + \beta \left(\left|v_l(T_0)\right| +\epsilon\right)^3\\
			\leq \, &\epsilon\left(2n+1+\frac{8n}{h^2}+2\alpha +8\beta
			\epsilon^2\right)<M.\end{align*}
		Thus, we have designed a control that drives the system to zero velocity at $T_1$.
		The system is now at a stationary flock at $t=T_1$.
		
		Finally, let 
		\[
		u_l(t)=-(\alpha-\beta|v_l|^2)v_l-\sum_{l' \sim l}\frac{x_{l'}(t)-x_l(t)}{h^2} + w \quad \text{ for $t \in [T_1,T_2]$,}
		\]
		where 
		\[
		T_2=T_1 + \frac{2\sqrt{\alpha}}{M\sqrt{\beta}} \quad \text{ and } \quad w=\frac{\bar{v}}{T_2-T_1}=\frac{M}{2}.
		\]
		We notice that the first terms of $u_l$ will cancel out the other terms on the right side of the second equation in \eqref{eq:sys} so that 
		\[
		\dot{v}_l(t)=w=\frac{M}{2}.
		\]
		Additionally, $v_l(T_1)=0$. 
		Thus, $v_l(t)$ is the same for all $l\in \overline \Lambda$ and
		\[
		v_l(t)=\frac{M(t-T_1)}{2} \quad \text{ for } t \in [T_1,T_2].
		\]
		Therefore, 
		\[
		\sum_{l' \sim l} \frac{x_{l'}(t)-x_l(t)}{h^2}=\sum_{l' \sim l} \frac{x_{l'}(T_1)-x_l(T_1)}{h^2} \quad \text{ for all $t \in [T_1, T_2]$.}
		\]
		Recall that this term has magnitude less than {$\epsilon(1+8n/h^2)< M/2$}. With this information, we can show that $|u_l|<M$. 
		As $0\leq v_l(t)\leq \bar v$ we can see using \eqref{fmax} that
		\begin{align*}
			0 \leq (\alpha - \beta |v_l(t)|^2)v_l(t) \leq  \frac{2}{3\sqrt{3}}\sqrt{\frac{\alpha ^3}{\beta}}<M.
		\end{align*}
		Thus,
		\[
		|u_l(t)|\leq \left|-(\alpha-\beta|v_l|^2)v_l+ \frac{M}{2}\right| +\left|-\sum_{l' \sim l}\frac{x_{l'}(t)-x_l(t)}{h^2} \right| \leq \frac{M}{2}+\epsilon{\left(1+\frac{8n}{h^2}\right)} <M.
		\]
		At $t=T_2$, we get
		\[
		v_l(T_2)=\bar v = \sqrt{\frac{\alpha}{\beta}},
		\]
		and
		\[
		\left|\sum_{l' \sim l} \frac{x_{l'}(T_2)-x_l(T_2)}{h^2}\right|=\left|\sum_{l' \sim l} \frac{x_{l'}(T_1)-x_l(T_1)}{h^2}\right| \leq \epsilon{\left(1+\frac{8n}{h^2}\right)}. 
		\]
		For $t\geq T_2$, we simply keep $v_l(t)=v_l(T_2)=\bar v$ for all $l\in \overline \Lambda$.
		Then, we get
		\[
		x_l(t) = x_l(T_2) + \bar v(t-T_2) \quad \text{ for all } t\geq T_2.
		\]
		The proof is complete with $T=T_2$ and $\bar x_{1l}=x_l(T_2)$ for $l\in \overline \Lambda$.
	\end{proof}

We note that the construction also gives an explicit bound on the time needed
after $T_0$. Since $|v_l(T_0)|<\epsilon$ for every
$l\in\overline{\Lambda}$, the definitions of $T_1$ and $T_2$ yield
\[
T_1-T_0
=\frac{1}{\epsilon}
\max_{l\in\overline{\Lambda}}|v_l(T_0)|<1,
\qquad
T_2-T_1=\frac{2\sqrt{\alpha}}{M\sqrt{\beta}}.
\]
Therefore, the moving flock is reached at time $T_2$, with
\begin{equation}\label{eq:postT0}
T_2-T_0<1+\frac{2\sqrt{\alpha}}{M\sqrt{\beta}}.
\end{equation}
\begin{correction}
    Note that no explicit bound on $T_0$ has been obtained.\end{correction}
	\begin{rem}
		It is important to note that, in the above proof, for $l\in \overline \Lambda$,
		{\[
		\left|\sum_{l' \sim l} \frac{x_{l'}(T_2)-x_l(T_2)}{h^2}\right|=\left|\sum_{l' \sim l} \frac{x_{l'}(T_1)-x_l(T_1)}{h^2}\right| \leq \frac{8n\epsilon}{h^2}+\epsilon,
		\]}
		where $\epsilon>0$ can be chosen as small as we want.
		This is an important point that we will need to use in the proof of Theorem \ref{thm:main2}.
	\end{rem}

	\section{Proofs of Theorems \ref{thm:main2}--\ref{thm:main3}}\label{sec:thm2-3}
	\subsection{Proof of Theorem \ref{thm:main2}}\label{sec:moving-proof}
	Our strategy to prove Theorem \ref{thm:main2} is as follows.
	By Theorem \ref{thm:main1}, we have driven the system to a moving flock at time $T_2$.
	By designing an explicit control, we drive the system to a moving flock where every particle is at the same position at $T_3>T_2$.
	From this intermediate moving flock where all particles are at the same position, we will finally design another explicit control driving the system to the desired flock at $T_4>T_3$.
	Along the way, we need to make sure that the changes of the positions of the particles and their corresponding velocities are sufficiently small so that the controls always satisfy the given bound $M$.
	
	\begin{proof}[Proof of Theorem \ref{thm:main2}]
		By the proof of Theorem \ref{thm:main1}, we know that given any initial conditions, we can drive the system to a moving flock $\varphi_1(l,t)=\bar \varphi_1(l) +\bar v(t{-T_2)}$ at time $t=T_2$ satisfying
		{\begin{equation} \label{epsbound} \left|\sum_{l' \sim l} \frac{\bar \varphi_1(l')-\bar \varphi_1(l)}{h^2}\right| \leq \epsilon\left(\frac{8n}{h^2}+1\right) \end{equation}} for all $l\in \overline \Lambda$.
		We will then adjust the control to drive the system to 
		\[
		\begin{cases}
			x_l(T_3)=x_j(T_3),\\
			v_l(T_3)=v_j(T_3)=\bar v,
		\end{cases}
		\]
		for all $l,j \in \overline \Lambda$.
		That is, the system is at a moving flock where every particle is at the same position at $T_3>T_2$.
		Let 
		\[
		\bar{y}=\frac{1}{D^n}\sum_{l\in \overline \Lambda} \bar{\varphi}_1(l),
		\]
		and let $T_3>T_2$ satisfy
		\begin{equation} \begin{split} \label{T2cond} &\frac{3\alpha |\bar{y} - \bar{\varphi}_1(l)|(T_3-T_2)^2+6.75\sqrt{\alpha \beta}|\bar{y}-\bar{\varphi}_1(l)|^2(T_3-T_2)}{(T_3-T_2)^3}\\
				&+\frac{3.375\beta |\bar{y} -\bar{\varphi}_1(l)|^3 +6|\bar{y}-\bar{\varphi}_1(l)|(T_3{-}T_2)}{(T_3-T_2)^3}<M-\epsilon{\left(\frac{8n}{h^2}+1\right)}\end{split}\end{equation}
		for all $l \in \overline\Lambda$.
\begin{rem}
    We note that our choice of $\bar{y}$ is reasonable, but another point could be chosen. In fact, if we choose $\bar{y}$ such that it minimizes $\max_{l \in \Lambda} |y-\bar\varphi_1(l)|$
then we can reduce $T_3$.
\end{rem}

		We then  set 
		\begin{equation*}\begin{split}
				u_l(t)=&-\left[1+\frac{2(t-T_2)^3}{(T_3-T_2)^3}-\frac{3(t-T_2)^2}{(T_3-T_2)^2}\right]\sum_{l'\sim l} \frac{\bar{\varphi}_1(l')-\bar{\varphi}_1(l)}{h^2} \\
				&+ \frac{12 \alpha}{(T_3-T_2)^3}(\bar{y}-\bar{\varphi}_1(l))(t-T_2)(T_3-t)\\
				& + \frac{108\sqrt{\alpha \beta}}{(T_3-T_2)^6}(\bar{y} - \bar{\varphi}_1(l))^2(t-T_2)^2(T_3-t)^2\\
				&+ \frac{216 \beta}{(T_3-T_2)^9}(\bar{y} - \bar{\varphi}_1(l))^3(t-T_2)^3(T_3-t)^3+\frac{6(\bar{y}-\bar{\varphi}_1(l))}{(T_3-T_2)^3}(T_3+T_2-2t)\end{split} \end{equation*}
		for $t\in [T_2,T_3]$.
		With this choice of control and the fact that $v(T_2)=\bar{v}$, we find
		\[v_l(t)=\bar{v} + \frac{6(\bar{y}-\bar{\varphi}_1(l))}{(T_3-T_2)^3}(t-T_2)(T_3-t)\] and
		\[x_l(t)=\bar{\varphi}_1(l)+(t-T_2)\bar{v}+ \frac{3(t-T_2)^2}{(T_3-T_2)^2}(\bar{y}-\bar{\varphi}_1(l))-\frac{2(t-T_2)^3}{(T_3-T_2)^3}(\bar{y}-\bar{\varphi}_1(l)).\]
		Since $x_l(T_3)=\bar{y} +(T_3-T_2) \bar{v}$, we deduce $x_l(T_3)=x_j(T_3)$ and $v_l(T_3)=v_j(T_3)=\bar v$ for all $l, j \in \overline{\Lambda}$. 
		However, we must show that $|u_l(t)|<M$. 
		Notice that $(t-T_2)(T_3-t)$ has a maximum value of $\frac{(T_3-T_2)^2}{4}$ and 
		\[
		1+\frac{2(t-T_2)^3}{(T_3-T_2)^3}-\frac{3(t-T_2)^2}{(T_3-T_2)^2}\leq 1 \quad \text{ for $t \in [T_2,T_3]$.}
		\]
		Combining these facts with \eqref{epsbound} and \eqref{T2cond} shows
		\begin{align*}
			&|u_l(t)|\\
			\leq\,& \left|\sum_{l' \sim l} \frac{\bar{\varphi}_1(l')-\bar{\varphi}_1(l)}{h^2}\right| + \frac{12 \alpha}{(T_3-T_2)^3}|\bar{y}-\bar{\varphi}_1(l)|\frac{(T_3-T_2)^2}{4} \\
			&\,+\frac{108 \sqrt{\alpha\beta}}{(T_3-T_2)^6}|\bar{y}-\bar{\varphi}_1(l)|^2\frac{(T_3-T_2)^4}{16} +\frac{216 \beta}{(T_3-T_2)^9}|\bar{y}-\bar{\varphi}_1(l)|^3 \frac{(T_3-T_2)^6}{64} \\
			&\,+ \frac{6{|}\bar{y} - \bar{\varphi}_1(l){|}}{(T_3-T_2)^3}(T_3{-}T_2)\\
			\leq\,& \epsilon{\left(\frac{8n}{h^2}+1\right)} + \frac{3\alpha |\bar{y} - \bar{\varphi}_1(l)|(T_3-T_2)^2+6.75\sqrt{\alpha\beta}|\bar{y}-\bar{\varphi}_1(l)|^2(T_3-T_2)}{(T_3-T_2)^3}\\
			&\,\,+\frac{3.375\beta |\bar{y} -\bar{\varphi}_1(l)|^3 +6|\bar{y}-\bar{\varphi}_1(l)|(T_3{-}T_2)}{(T_3-T_2)^3}<M.
		\end{align*}
		From the intermediate moving flock where all particles are at the same position, we will finally design the control for $t \in [T_3, T_4]$ that drives the system to the desired flock
		\[
		x_l(t)= \begin{correction} \bar{\varphi}_2(l) +\bar{v}t\end{correction} \quad \text{ for } l \in \overline \Lambda, t \geq T_4.
		\]
		By the assumptions, $|-\Delta \bar{\varphi}_2|<M$. 
        \begin{correction} Let $d_l=\bar{\varphi}_2(l)+T_3\bar{v}-x_l(T_3)$.\end{correction}
		Let $T_4>T_3$ satisfy
		\begin{equation} 
			\begin{split} \label{T3cond} 
				&\frac{3\alpha \begin{correction}|d_l| \end{correction}(T_4-T_3)^2+6.75\sqrt{\alpha\beta}\begin{correction}|d_l| \end{correction}^2(T_4-T_3)+3.375\beta \begin{correction}|d_l| \end{correction}^3 +6\begin{correction}|d_l| \end{correction}(T_4{-}T_3)}{(T_4-T_3)^3}\\
                &<M-|\Delta \bar{\varphi}_2(l)|
			\end{split}
		\end{equation}
		for all $l \in \overline \Lambda$.
		
		\smallskip
		
		Now from $T_3$ to $T_4$, letting
		\begin{equation*} 
			\begin{split} \label{nextu}
				u_l(t)=&\left[-\frac{3(t-T_3)^2}{(T_4-T_3)^2} + \frac{2(t-T_3)^3}{(T_4-T_3)^3}\right] \sum_{l' \sim l} \frac{\bar{\varphi}_2(l')-\bar{\varphi}_2(l)}{h^2}\\
				&+ \frac{12 \alpha}{(T_4-T_3)^3}\begin{correction}d_l \end{correction}(t-T_3)(T_4-t) + \frac{108\sqrt{\alpha \beta}}{(T_4-T_3)^6}\begin{correction}d_l \end{correction}^2(t-T_3)^2(T_4-t)^2\\
				&+ \frac{216 \beta}{(T_4-T_3)^9}\begin{correction}d_l\end{correction}^3(t-T_3)^3(T_4-t)^3+\frac{6\begin{correction}d_l \end{correction}}{(T_4-T_3)^3}(T_4+T_3-2t),
			\end{split} 
		\end{equation*}
		we arrive at
		\[v_l(t)=\bar{v} + \frac{6\begin{correction}d_l\end{correction}}{(T_4-T_3)^3}(t-T_3)(T_4-t)\]
		and
		\[x_l(t)=x_l(T_3) + (t-T_3)\bar{v} + \frac{3(t-T_3)^2}{(T_4-T_3)^2}\begin{correction}d_l \end{correction}-\frac{2(t-T_3)^3}{(T_4-T_3)^3}\begin{correction}d_l \end{correction}.\]
		Then,
		\[
		\begin{cases}
			x_l(T_4)=\begin{correction}T_4\bar{v}\end{correction}+\bar{\varphi}_2(l),\\
			v_l(T_4)=\bar{v}.
		\end{cases}
		\]
		Once again, we would like to ensure $|u_l(t)|<M$. 
		For $l\in \overline \Lambda$, we  estimate
		\begin{align*}
			|u_l(t)|
			\leq\,  &|\Delta \bar{\varphi}_2(l)| + \frac{3\alpha \begin{correction}|d_l| \end{correction}(T_4-T_3)^2+6.75\sqrt{\alpha\beta}\begin{correction}|d_l| \end{correction}^2(T_4-T_3)+3.375\beta \begin{correction}|d_l| \end{correction}^3}{(T_4-T_3)^3}\\
            &+\frac{6\begin{correction}|d_l| \end{correction}(T_4{-}T_3)}{(T_4-T_3)^3}<M.
		\end{align*}
		Finally, for $t>T_4$, we will let
		\[u_l(t)=-\sum_{l' \sim l} \frac{x_{l'}(t)-x_l(t)}{h^2},\]
		so $v_l(t)=\bar{v}$, and $x_l(t)=\begin{correction}\bar{v}t + \bar{\varphi}_2(l)\end{correction}$ for $t\geq T_4$.
	\end{proof}

\subsection{Proof of Theorem \ref{thm:main3}}\label{sec:stationary-proof}
\begin{proof}[Proof of Theorem \ref{thm:main3}]
For a lattice vector $q$, we define
$\|q\|_\infty=\max_{l\in\overline\Lambda}|q_l|$.
We now use the cubic interpolation polynomial
\[
p(s)=3s^2-2s^3,\qquad 0\leq s\leq1,
\]
which satisfies
\[
p(0)=0,\quad p(1)=1,\quad p'(0)=p'(1)=0,
\quad 0\leq p\leq1,\quad |p'|\leq\frac32,\quad |p''|\leq6.
\]
We now use the first two phases of Section~\ref{sec:thm1} and stop at $T_1$. 
We set $A_l=x_l(T_1)$.

Then, we obtain
\begin{equation}\label{epsbound1}
v_l(T_1)=0,\qquad \|\Delta x_l(T_1)\|_\infty\leq E_h<M/2,
\end{equation}
where $E_h:=\epsilon(1+8n/h^2)$.
We let
\[
c=\frac1{D^n}\sum_{l\in\overline{\Lambda}} x_l(T_1),\qquad d_l=c-x_l(T_1),
\qquad R_5=\|d\|_\infty,\qquad \mu_5=M-E_h>0.
\]
For $R\geq0$, we can define the bound
\begin{equation}\label{eq:Qstationary}
Q_{\rm s}(\tau,R)=\frac{3\alpha R}{2\tau}
                  +\frac{6R}{\tau^2}
                  +\frac{27\beta R^3}{8\tau^3}.
\end{equation}
We can compute an explicit time bound for  $\mu>0$  
\begin{equation}\label{eq:stationary-duration}
\tau_{\rm s}(R,\mu)=\max\left\{1,\frac{6\alpha R}{\mu},
                    \sqrt{\frac{24R}{\mu}},
                    \left(\frac{27\beta R^3}{2\mu}\right)^{1/3}\right\}.
\end{equation}
We note that each summand in \eqref{eq:Qstationary} is then at most $\mu/4$, whence
$Q_{\rm s}(\tau_{\rm s}(R,\mu),R)\leq3\mu/4<\mu$.

We now choose $\tau_5=\tau_{\rm s}(R_5,\mu_5)$ and $T_5=T_1+\tau_5$.
For $T_1\leq t\leq T_5$, we then set $s=(t-T_1)/\tau_5$ and prescribe
\begin{align}
x_l^*(t)&=A_l+p(s)d_l,\qquad
v_l^*(t)=\frac{p'(s)}{\tau_5}d_l,\label{eq:stationary-path1}\\
u_l(t)&=\frac{p''(s)}{\tau_5^2}d_l
        -(1-p(s))\Delta x_l(T_1)
        -\alpha\frac{p'(s)}{\tau_5}d_l
        +\beta\left(\frac{p'(s)}{\tau_5}d_l\right)^3.
\label{eq:stationary-control1}
\end{align}
Here $\Delta x^*=(1-p(s))\Delta x_l(T_1)$ as $c$ is common to all sites.
Thus, \eqref{eq:stationary-control1} is exactly
$u_l=\dot v_l^*-\Delta x_l^*-f(v_l^*)$, and the prescribed state solves the controlled system. Uniformly in $l$ and $t$, we have the bound
\begin{equation}\label{eq:stationary-bound1}
|u_l(t)|\leq E_h+Q_{\rm s}(\tau_5,R_5)
           \leq E_h+\tfrac34\mu_5<M,
\end{equation}
which implies $x_l(T_5)=c$ and $v_l(T_5)=0$.

For the second stage, \begin{correction} we choose 
\[
e_l=\bar\varphi_3(l)-c,
\qquad R_6=\|e\|_\infty,\qquad \mu_6=M-\|\Delta \bar\varphi_3(l)\|_\infty>0.
\]
\end{correction}
We now set
$\tau_6=\tau_{\rm s}(R_6,\mu_6)$ and $T_6=T_5+\tau_6$.
For $T_5\leq t\leq T_6$, we can put $s=(t-T_5)/\tau_6$ and prescribe
\begin{align}
x_l^*(t)&=c+p(s)e_l,\qquad
v_l^*(t)=\frac{p'(s)}{\tau_6}e_l,\label{eq:stationary-path2}\\
u_l(t)&=\frac{p''(s)}{\tau_6^2}e_l
        -p(s)\Delta \bar\varphi_3(l)
        -\alpha\frac{p'(s)}{\tau_6}e_l
        +\beta\left(\frac{p'(s)}{\tau_6}e_l\right)^3.
\label{eq:stationary-control2}
\end{align}
We obtain the uniform control estimate 
\begin{equation}\label{eq:stationary-bound2}
|u_l(t)|\leq\|\Delta \bar\varphi_3(l)\|_\infty+Q_{\rm s}(\tau_6,R_6)
           \leq\|\Delta \bar\varphi_3(l)\|_\infty+\tfrac34\mu_6<M.
\end{equation}

At $T_6$, $x_l(T_6)=\begin{correction}\bar\varphi_3(l)\end{correction}$ and $v_l(T_6)=0$. For $t\geq T_6$, the constant control $u_l=-\Delta B_l$ maintains this state and satisfies $|u_l|<M$.  The resulting bounded, piecewise continuous control belongs to $\cA_M$, and the theorem follows with $T=T_6$.
\end{proof}

\begin{rem}\label{rem:stability}
Let us  fix a terminal profile $\varphi$ and retain the prescribed force $u_l=-\Delta \varphi(l)$. We then set 
\[
x_l(t)=\varphi(l)+\bar v t+y_l(t),\qquad v_l(t)=\bar v+z_l(t),
\]
where $\bar v\in\{0,\pm\sqrt{\alpha/\beta}\}$. The exact perturbation equations are
\begin{equation}\label{eq:perturbation}
\dot y_l=z_l,\qquad
\dot z_l=\Delta y_l+(\alpha-3\beta\bar v^2)z_l
                      -3\beta\bar v z_l^2-\beta z_l^3.
\end{equation}

The finite periodic lattice is connected, and thus $-\Delta$ is a symmetric and nonnegative operator, with a one-dimensional kernel consisting of constant vectors. For an eigenmode of $-\Delta$ with eigenvalue $\lambda\geq0$, the linearized amplitude $\eta$ satisfies
\[
\ddot\eta-(\alpha-3\beta\bar v^2)\dot\eta+\lambda\eta=0.
\]
For a moving flock, we find $$\ddot\eta+2\alpha\dot\eta+\lambda\eta=0.$$ If $\lambda>0$, both characteristic roots have negative real parts. If $\lambda=0$, the roots are $0$ and $-2\alpha$. Thus, a moving flock is linearly asymptotically stable modulo this neutral translation mode.

For a stationary flock, the homogeneous mode satisfies
$$\ddot\eta-\alpha\dot\eta=0$$ and has the positive characteristic root $\alpha$. Nonlinear instability can also be seen directly: for spatially homogeneous $y,z$, \eqref{eq:perturbation} reduces to
$$\dot z=\alpha z-\beta z^3.$$ Every sufficiently small positive initial velocity increases toward $\sqrt{\alpha/\beta}$ rather than returning to zero. Hence, a stationary flock is  unstable under the fixed target force. 
\end{rem}

\section{Connection to optimal control and Hamilton--Jacobi equations}\label{sec:connect}
    The goal of the control strategy above is to drive the system to a flock in a finite amount of time. Although the strategy strives to do this efficiently, it is by no means optimal. A related problem is to find the minimal time it takes to start from a given state and reach a flock with a control that has magnitude {at most} $M$. There is a rich theory behind this kind of optimal control problem. Here, we connect this optimal problem to a Hamilton-Jacobi problem (see \cite{Tran} for an overview on the topic).
    
	{The first step is to rewrite the coupled lattice \eqref{eq:sys} as a finite-dimensional controlled ODE. } 
	Thus, we write $X_1=(x_l)_{l \in \ol{\Lambda}}$ and $X_2=(v_l)_{l\in \ol{\Lambda}}$ as vectors in $\R^{D^n}$.
	Denote by $X=(X_1,X_2)\in \R^{2 D^n}$.
	Then, \eqref{eq:sys} can be written as
	\[
	\dot X = \tilde b(X,W)= b(X)+ (0,W),
	\]
	with the given initial condition
	\[
	X(0)=\left( (w_0(l))_{l \in \ol{\Lambda}},(w_1(l))_{l \in \ol{\Lambda}}\right).
	\]
	Here, the vector field $b$ is
	\[
	b(X) = \left(X_2,  \left(\sum_{l' \sim l} \frac{(X_1)_{l'} - (X_1)_l}{h^2} + (\alpha - \beta |(X_2)_l|^2) (X_2)_l\right)_{l\in \ol{\Lambda}} \right).
	\]
	In this framework, $W$ is the control, and the admissible set is
	\[
	\ol\cA_M=\left\{W\in \R^{D^n}\,:\,|W_l| =|u_l| \leq M \qquad \text{ for } l\in \ol{\Lambda} \right\}.
	\]
	Basically, we can think of $X$ as the master state since it contains all information, and {the ODE for $X$ retains all lattice interactions in $b(X)$}.
	
	\smallskip
	
	Without loss of generality, we are only concerned with moving flocks here.
	For a given state $X$, we are interested in the optimal time it takes to reach a moving flock, that is, 
	\begin{multline*}
		\Gamma = \left\{X\;:\, \left|\sum_{l' \sim l} \frac{(X_1)_{l'} - (X_1)_l}{h^2}\right| \leq M, (X_2)_l = -\sqrt{\alpha/\beta} \quad \text{ for } l\in \ol{\Lambda}\right\}\\ 
		\bigcup \left\{X\;:\, \left|\sum_{l' \sim l} \frac{(X_1)_{l'} - (X_1)_l}{h^2}\right| \leq M, (X_2)_l = \sqrt{\alpha/\beta} \quad \text{ for } l\in \ol{\Lambda}\right\}.
	\end{multline*}
	Let $U(X)$ be the minimal time it takes to start from state $X$ to reach $\Gamma$.
	Then, for $X=(X_1,X_2)\in \R^{2 D^n}$ and $P=(P_1,P_2)\in \R^{2 D^n}$, the corresponding Hamiltonian is
	\[
	H(X,P) = \sup_{W\in \ol \cA_M} \left( -\tilde b(X,W)\cdot P -1\right)=-b(X)\cdot P + M|P_2|_1-1.
	\]
	Here, $|\cdot|_1$ is the usual $1$-norm in $\R^{D^n}$.
	{The formal dynamic programming equation associated with this minimum-time problem is the following HJ-type equation: )}
	\begin{equation}\label{eq:eikonal}
		\begin{cases}
			-b(X)\cdot D_X U(X) + M |D_{X_2} U(X)|_1 =1 \qquad &\text{ for } X \in \R^{2 D^n}{\setminus\Gamma},\\
			U(X)=0 \qquad &\text{ for } X \in \Gamma.
		\end{cases}
	\end{equation}
	{For viscosity solutions and the minimal time functions, see \cite[Chapter 2]{Tran}.}
	We note that this is a very complicated PDE posed in high dimensions, as the vector field $b$ contains much information about the dynamics.
	The value function $U(X)$ gives us the minimal time to reach a moving flock starting from $X\in \R^{2 D^n}$, but does not give us directly an optimal strategy to reach a moving flock.
    By Theorem \ref{thm:main1}, we have that $U$ is finite.
    As the Hamiltonian $H$ here is noncoercive in $P$, the regularity of the {minimum-time function} $U$ is not yet known.
    {In particular, continuity and boundary regularity of $U$ remain open. }
	
	The advantage of Theorems \ref{thm:main1}--\ref{thm:main2} is that they provide more or less explicit strategies to reach a moving flock in finite time. Nevertheless, it is important to emphasize that an optimal control framework is hidden behind this problem.
    
    Equation \eqref{eq:eikonal} is interesting in its own right and merits further study.
    Specifically, the regularity of the {minimum-time function} $U$, the properties of the optimal trajectories (strategies), and possible approximations are of particular interest.

\section{Numerical illustrations}\label{sec:numerics}
The computation below provides examples of our control strategy illustrated above. \begin{correction}These examples were run using Python 3.13.3, NumPy 2.5.4, SciPy 1.18.1, Matplotlib 3.11.2, and pandas 3.1.0. The script is available publicly on Zenodo \cite{STTCode}.\end{correction}

We use a periodic one-dimensional chain with $D=12$, $h=1/12$, $\alpha=\beta=1$, and sites $j=0,\ldots,11$. The initial data and target {at the lattice site $l=jh$} are
\begin{align*}
x_l(0)&=0.04\cos(2\pi {l})+0.015\sin(4\pi {l}),\\
v_l(0)&=-1/\sqrt3+0.12\sin(2\pi 
{l}),\\
 \varphi(l)&=0.008\cos(2\pi {l}).
\end{align*}
We run the procedure with multiple values for $M$. In each case, we choose $\gamma=2\max\{1,1/M\}$ and take $\epsilon=\frac{M}{26+16/h^2}$. For the choices of $M$ and $h$ being  used, this ensures \eqref{eps} is satisfied. The first phase is integrated using the explicit eighth-order DOP853 solver, with relative tolerance $10^{-10}$, absolute tolerance $10^{-12}$, and maximum step $0.05$. The {initial Lyapunov-feedback phase} described at the beginning of Section \ref{sec:thm1} is carried out (see \eqref{ui1}) until the Lyapunov functional satisfies
\[
V=\frac{\epsilon^2}{8\max\{1,4/h^2\}}=\frac{\epsilon^2h^2}{32}.
\]

Our next phase is completed until $\max_{l \in \overline{\Lambda}} T_{1,l}$, and then the system is driven to a flock, {with positions $\phi_1(l):=x_l(T_2)$ where} $T_2=T_1 + 2/M$. {Following   the proof of Theorem~\ref{thm:main2}, we also set $$\bar y=D^{-1}\sum_{l\in\overline\Lambda}\phi_1(l).$$} $T_3$ and $T_4$ must satisfy \eqref{T2cond} and \eqref{T3cond}, respectively. {The main example uses $M=1$, while the threshold study below uses $M\leq 2F_*$. In these examples, we set}
\[T_3=\max_{l \in \overline{\Lambda}}\left(\max\left\{\frac{6|\bar{y}-{\phi_1(l)}|}{M-\epsilon(8n/h^2+1)}, \sqrt{\frac{24|\overline{y}-{\phi_1(l)}|}{M-\epsilon(8n/h^2+1)}}\right\}\right)+T_2\]
and
\[T_4=\max_{l \in \overline{\Lambda}}\left(\max\left\{\frac{6|\varphi(l)-x_l(T_3)|}{M-\max_{l \in \overline{\Lambda}}|\Delta \varphi(l)|}, \sqrt{\frac{24|\varphi(l)-x_l(T_3)|}{M-\max_{l \in \overline{\Lambda}}|\Delta \varphi(l)|}}\right\}\right)+T_3.\]

The following figure shows the details of the trajectories of $x_l$, $v_l$, and $u_l$ when $M=1$.

\begin{figure}[htbp]
\centering
\includegraphics[width=\textwidth]{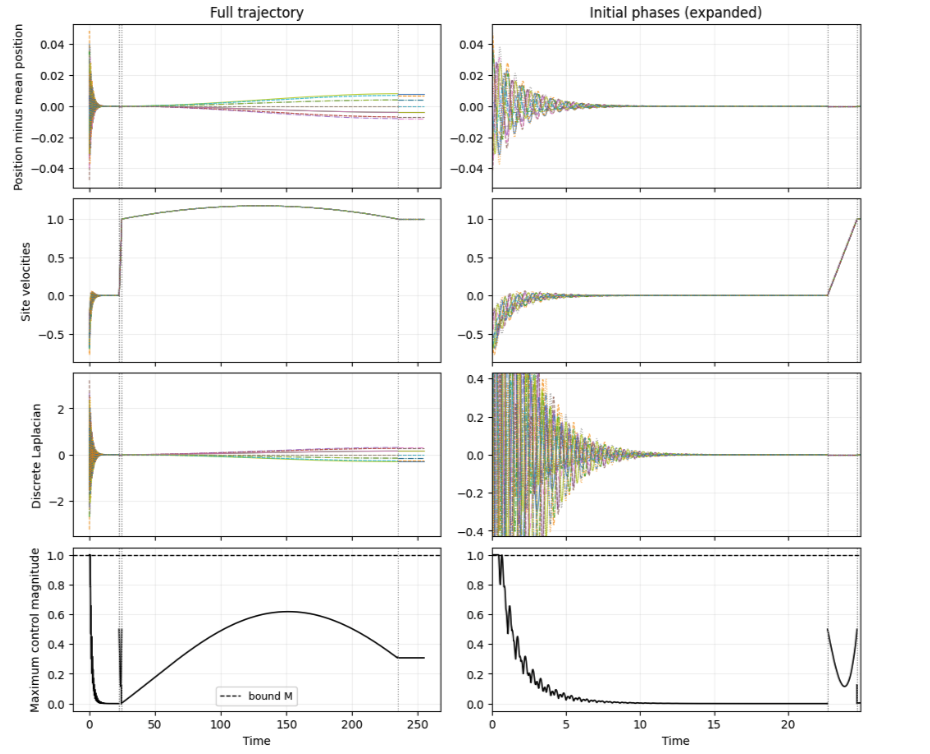}
\caption{\begin{correction}Moving-target construction with $M=1$. The left panels show trajectories over all time and the right panels show trajectories in the initial stages of the construction. The top row corresponds to the positions minus the mean of all sites. The second row shows the velocities at each site. The third row shows the discrete Laplacian values and the fourth row shows $\max_l|u_l|$ and the bound $M$.\end{correction} Vertical lines mark $T_0$, $T_1$, $T_2$, {$T_3$, and $T_4$; nearby switches may appear coincident at this scale}.}\label{fig:phases}
\end{figure}

For $M=1$, the switching times are chosen to be $T_0=22.649432$, $T_1=22.659214$, $T_2=24.659214$, $T_3=24.662854$, $T_4=28.750748$. The last two stage durations can be computed  $T_3-T_2=0.003640$ and $T_4-T_3=4.087894$. The maximum sampled control magnitude is $1$. After $T_0$ it is $0.5$. 
For a threshold study, we set $M=F_*(1+r)$ with
\[
r\in\{0.002,0.005,0.01,0.02,0.05,0.1,0.2,0.5,1\}
\]
and retain the same initial data and target. The target remains strictly admissible at every tested $M$.   The first phase could be stopped at a different time, which would change the results. We note that the construction described in this paper brings the configuration to a \begin{correction}prescribed flock\end{correction}.
\begin{figure}[htbp]
\centering
\includegraphics[width=.92\textwidth]{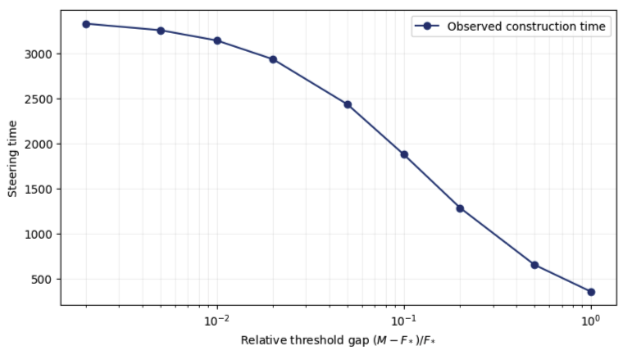}
\caption{Construction time versus the relative distance above the Rayleigh-force threshold. A logarithmic x-axis is used.}\label{fig:time}
\end{figure}

\clearpage

\section{Conclusion}\label{sec:conclusion}
We consider a finite periodic lattice with fixed mesh size and an independent scalar control at each site. Under Assumption~1, our construction steers every initial state to a moving flock in finite time. Under the hypotheses of Theorems~\ref{thm:main2}--\ref{thm:main3}, it also reaches a \begin{correction} prescribed moving or stationary profile.\end{correction} One of our main analytical contributions is an explicit, phase-by-phase construction that enforces the pointwise control bound. 

The construction provides explicit sufficient controls without asserting time optimality. The initial feedback phase brings the system into a prescribed small neighborhood in finite time, but we do not obtain a general explicit bound on $T_0$.  

The trajectory-planning phases require exact knowledge of the model parameters and exact initialization at each switching time. Remark~\ref{rem:stability} presents the coupled perturbation equations and shows that, under the fixed target force, moving flocks are linearly stable modulo common spatial translations, whereas stationary flocks are unstable. \begin{correction} We leave the study of nonperiodic lattices or sparse controls for future study. We also note that this work does not allow $h\to0$. However, for the application studied in this paper, $h$ is fixed and should not be sent to zero.
\end{correction}

\section{Statement of AI use}\label{sec:ai}
\begin{correction}
All proofs and work in Sections 1-4 were completed without the use of AI. 

ChatGPT (OpenAI) was used to assist with drafting and debugging the Python implementation and preparing numerical illustrations. The authors reviewed the code and take full responsibility for the content. The computations reported in Section~\ref{sec:numerics} were executed in Google Colab using NumPy, SciPy, and Matplotlib. 
The code is available at \\
\url{https://github.com/slstrikwerda/Controlling-discrete-semilinear-wave-equations-toward-flocks}.
\end{correction}

	\end{document}